\documentclass[12pt]{article}
\usepackage{amsmath}
\usepackage{amssymb}
\renewcommand\smallskip{\vskip\smallskipamount}
\renewcommand\medskip{\vskip\medskipamount}
\renewcommand\bigskip{\vskip\bigskipamount}

\begin{document}

\footnotetext{The author is partially supported by NSF Grant
DMS-0707086 and a Sloan Research Fellowship.}

\begin{center}
\begin{large}
\textbf{The Hoop Conjecture in Spherically Symmetric Spacetimes}
\end{large}

\bigskip\smallskip

MARCUS A. KHURI

\bigskip

\end{center}

\begin{abstract}
We give general sufficient conditions for the existence of trapped
surfaces due to concentration of matter in spherically symmetric
initial data sets satisfying the dominant energy condition.  These
results are novel in that they apply and are meaningful for
arbitrary spacelike slices, that is they do not require any
auxiliary assumptions such as maximality, time-symmetry, or
special extrinsic foliations, and most importantly they can easily
be generalized to the nonspherical case once an existence theory
for a modified version of the Jang equation is developed.
Moreover, our methods also yield positivity and monotonicity
properties of the Misner-Sharp energy.
\end{abstract}

\begin{center}
\textbf{1. Introduction}
\end{center}\setcounter{equation}{0}
\setcounter{section}{1}

  The Hoop Conjecture [1] concerns the folklore belief that if
enough matter and/or gravitational energy are present in a small
enough region (small in all three spatial dimensions), then the
system must collapse to a black hole.  This belief is often
realized by establishing a statement of the following form.  Let
$\Omega$ be a compact spacelike hypersurface satisfying an
appropriate energy condition in a spacetime $\mathcal{M}$.  There
exists a universal constant $C>0$ such that if
$\mathrm{Mass}(\Omega)>C\cdot\mathrm{Size}(\Omega)$, then $\Omega$
must contain a closed trapped surface.  Of course finding the
correct notions of $\mathrm{Mass}(\Omega)$ and
$\mathrm{Size}(\Omega)$ is one of the primary difficulties with
this conjecture.  The conclusion of the above statement guarantees
that the spacetime $\mathcal{M}$ contains a singularity (or more
precisely is null geodesically incomplete) by the Hawking-Penrose
Singularity Theorems [2], and assuming cosmic censorship must
therefore contain a black hole. It should also be pointed out that
modulo certain technical restrictions, trapped surfaces lead to
gravitational confinement according to Israel's result [3].
Therefore, in asymptotically flat spacetimes the existence of a
trapped surface almost certainly implies the existence of a black
hole.\par
  There have been many results realizing a version of the hoop
conjecture in this spirit.  Most notable are those of O'Murchadha,
Malec, and others [4,5,6,7,8,9,10,\linebreak 11,12], which address
concentration of matter in spherical symmetry and give necessary
and sufficient conditions in some instances, but impose auxiliary
conditions on the spacelike slice such as the condition of
maximality, time-symmetry, or that it arises from an extrinsic
foliation. On the other hand, there are the very important results
of Schoen and Yau [13], [14] which also address the issue of
concentration of matter, but without extra assumptions. While
their results are very impressive in that they do not require
spherical symmetry, they suffer from the opposite problem in that
they are not meaningful for slices with small extrinsic curvature,
in particular for maximal or time-symmetric slices. There have
been relatively fewer results on the concentration of pure
gravitational radiation, see [15] and [9].\par
  In this paper we also address the topic of concentration of
matter.  Our goal is to establish sufficient conditions for the
existence of trapped surfaces in spherically symmetric initial
data, which apply and are meaningful both in the maximal and
general cases.  Our methods are quite general in that they can
easily be generalized to the nonspherically symmetric case once an
appropriate existence theory (analogous to that developed by
Schoen and Yau in [16]) for a modified version of the Jang
equation has been established.  Moreover our techniques yield
positivity and monotonicity properties for the Misner-Sharp
energy, as a natural and interesting byproduct.\par
  An initial data set for the gravitational field consists of a
3-manifold $M$ on which is defined a positive definite metric $g$
and symmetric 2-tensor $k$ representing the extrinsic curvature.
The metric and extrinsic curvature must of course satisfy the
constraint equations:
\begin{eqnarray*}
16\pi\mu&=&R+(\mathrm{Tr}_{g}k)^{2}-|k|^{2},\\
8\pi J_{i}&=&\nabla^{j}(k_{ij}-(\mathrm{Tr}_{g}k)g_{ij}),
\end{eqnarray*}
where $R$ denotes scalar curvature and $\mu$, $J$ are respectively
the energy and momentum densities for the matter fields.  If the
initial data are spherically symmetric then we may take
$M\simeq\mathbb{R}^{3}$ and write
\begin{eqnarray*}
g&=&g_{11}(r)dr^{2}+\rho^{2}(r)d\chi^{2},\\
k^{ij}&=&n^{i}n^{j}k_{a}(r)+(g^{ij}-n^{i}n^{j})k_{b}(r),
\end{eqnarray*}
where
\begin{equation*}
n=n^{1}\partial_{r}+n^{2}\partial_{\psi^{2}}+n^{3}\partial_{\psi^{3}}=
\sqrt{g^{11}}\partial_{r}
\end{equation*}
is the unit normal to spheres $S_{r}$ centered at the origin (the
ball enclosed by $S_{r}$ will be denoted $B_{r}$), and
\begin{equation*}
d\chi^{2}=(d\psi^{2})^{2}+\sin^{2}\psi^{2}(d\psi^{3})^{2}
\end{equation*}
is the round metric on $\mathbb{S}^{2}$.  We also assume that the
metric is regular at the origin, so that $\rho(0)=0$,
$\rho_{,r}(0)=1$, and $g_{11}(0)=1$.  The sphere $S_{r}$ is future
(past) trapped if the family of outgoing future (past) directed
null geodesics, orthogonal to $S_{r}$, is converging at each
point. This is equivalent to the following condition satisfied by
the the null expansions:
\begin{eqnarray*}
\theta_{+}&:=&H_{S_{r}}+\mathrm{Tr}_{S_{r}}k<0\text{ }\text{
}\text{
}\text{ (future trapped)},\\
\theta_{-}&:=&H_{S_{r}}-\mathrm{Tr}_{S_{r}}k<0\text{ }\text{
}\text{ }\text{ (past trapped)},
\end{eqnarray*}
where $H_{S_{r}}=n_{;i}^{i}$ denotes the mean curvature and
$\mathrm{Tr}_{S_{r}}k=(g_{ij}-n_{i}n_{j})k^{ij}$.  The outer
boundary of a region in $M$ which contains future (past) trapped
surfaces is called a future (past) apparent horizon, and satisfies
$\theta_{+}=0$ ($\theta_{-}=0$).  Our main result is the
following\medskip

\textbf{Theorem 1.}  \textit{Let $(M,g,k)$ be a spherically
symmetric initial data set satisfying the dominant energy
condition $\mu\geq |J|$.  If}
\begin{equation}
\min_{B_{r}}(\mu\mp J(n))+\frac{3}{32\pi}\theta_{+}\theta_{-}(r)>
\frac{3}{2}\frac{\mathrm{Rad}(B_{r})}{\mathrm{Vol}(B_{r})}
\end{equation}
\textit{where the radius and volume are given by}
\begin{equation*}
\mathrm{Rad}(B_{r})=\int_{0}^{r}\sqrt{g_{11}},\text{ }\text{
}\text{ }\text{
}\mathrm{Vol}(B_{r})=4\pi\int_{0}^{r}\sqrt{g_{11}}\rho^{2},
\end{equation*}
\textit{then $B_{r}$ contains a future (past) trapped
surface.}\medskip

The first term on the left-hand side of (1.1) shows the
intuitively obvious fact, that formation of trapped surfaces
depends not only on matter concentration but also on the direction
that the matter is flowing.  Namely, inward flowing material
hastens (delays) the formation of future (past) trapped surfaces,
whereas outgoing material delays (hastens) formation.  More
interesting is the second term on the left-hand side, which
indicates that the bending of light rays at the boundary of
$B_{r}$, can by itself cause surfaces to be trapped on the
interior.  This phenomenon was first observed by Yau [14]. However
as we have pointed out, the result of [14] as well as the earlier
version of Schoen and Yau [13] are not meaningful when
$\mathrm{Tr}_{g}k$ is small.  To be more precise let us recall
their result, which states that if
\begin{equation*}
\min_{B}(\mu-|J|)>\frac{3\pi^{2}}{2}\frac{1}{\mathrm{Rad}_{SY}(B)}
\end{equation*}
then $B$ contains a trapped surface.  Here spherical symmetry is
not assumed and $\mathrm{Rad}_{SY}(B)$ is the square of a
``homotopy radius".  Thus their result requires matter density to
be large on a ``large region".  However our basic intuition
suggests that this is not the ideal situation which results in
collapse, that is, as the hoop conjecture asserts we would rather
like to show that if matter density is large on a ``small region"
then trapped surfaces exist.  So it is not surprising that their
result is only meaningful for a fairly restricted class of initial
data (as pointed out by Bizon, Malec, and O'Murchadha [4]).  They
show (Theorem 1 of their paper) that one cannot have a large set
with large positive scalar curvature. Since the matter density is
related to the scalar curvature via the Hamiltonian constraint,
the only way we can have a large matter density and small scalar
curvature (which is required from their Theorem 1) is that the
trace of the extrinsic curvature $\mathrm{Tr}_{g}k$ is large; in
fact the trace must be not only large but significantly larger
than $|k|$.  This means that it may be difficult to find data
which satisfy their condition, and in particular, their result can
say nothing about the time-symmetric or maximal cases. On the
other hand, our result compares nicely with that of Malec and
O'Murchadha [17] who showed that under the assumption of spherical
symmetry and maximality ($\mathrm{Tr}_{g}k=0$),
\begin{equation*}
4\pi\int_{0}^{r}(\mu\mp J(n))\rho^{2}>\mathrm{Rad}(B_{r})
\end{equation*}
implies that $B_{r}$ contains a future (past) trapped surface.
Unfortunately, it is difficult to see how their arguments might
generalize to the nonspherically symmetric case.

\begin{center}
\textbf{2. The Generalized Jang Equation}
\end{center}\setcounter{equation}{0}
\setcounter{section}{2}

  Our methods are based in large part on the generalized Jang
equation [18], which we now explain.  Many of the difficult issues
and questions involving initial data are easier to express and
solve if it happens that the scalar curvature of the given metric
$g$ is nonnegative, $R\geq 0$.  Unfortunately there is no
guarantee that this will be the case for an arbitrary set of
initial data, except under the added assumptions of maximality and
the dominant energy condition.  It is for this reason that Jang
[19] introduced the quasilinear elliptic equation for a scalar $f$
depending on $g$ and $k$, which bears his name:
\begin{equation}
H_{\Sigma}-\mathrm{Tr}_{\Sigma}K=0,
\end{equation}
where $\Sigma$ denotes the graph $t=f(x)$ inside the product
manifold $(M\times\mathbb{R}, g+dt^{2})$, $H_{\Sigma}$ is the mean
curvature, and $K$ is a trivially extended version of $k$
(extended to all of $M\times\mathbb{R}$).  That is, he showed that
if $f$ solves (2.1) then the scalar curvature of the related
metric $\overline{g}=g+df^{2}$ (this is the induced metric on the
graph $\Sigma$) has nice positivity properties.  In fact, Schoen
and Yau [16] successfully employed the Jang equation in their
solution of the positive energy conjecture, to reduce the general
case to the case of time-symmetry.  Moreover they developed a full
existence theory for this equation, and showed that regular
solutions always exist if the initial data do not contain apparent
horizons.  The converse statement, that if a regular solution does
not exist then the data must contain an apparent horizon,
naturally led to their result [13] concerning the hoop
conjecture.\par
  These successful applications of the Jang equation led many to
suggest that it could also be used to study the Penrose
Inequality.  However as pointed out by Malec and O'Murchadha [10],
serious and immediate difficulties arise when attempting such an
application.  These difficulties motivated the author together
with H. Bray [18] to propose a modified version of the Jang
equation, specifically designed for the Penrose Inequality.  This
generalized Jang equation has the same geometric structure as that
of (2.1), however the mean curvature of the graph $\Sigma$ is now
calculated inside the warped product manifold
$(M\times\mathbb{R},g+\phi^{2} dt^{2})$ where $\phi$ is a
nonnegative scalar, and the extended tensor $K$ is now a
nontrivial extension of $k$ (see [20]). An important feature of
the generalized Jang equation, like the original, is that it
yields nice positivity properties for the scalar curvature of the
induced metric on $\Sigma$.  More precisely, if $\overline{R}$
denotes the scalar curvature of $\overline{g}=g+\phi^{2} df^{2}$
then we find ([18]) that
\begin{equation}
\overline{R}=16\pi(\mu-J(w))+|h-K|_{\Sigma}|_{\overline{g}}^{2}+2|q|_{\overline{g}}^{2}
-2\phi^{-1}\mathrm{div}_{\overline{g}}(\phi q),
\end{equation}
where $h$ is the second fundamental form of $\Sigma$, and the
1-forms $w$ and $q$ are given by
\begin{equation*}
w_{i}=\frac{f_{,i}}{\sqrt{\phi^{-2}+|\nabla_{g}f|^{2}}},\text{
}\text{ }\text{ }\text{ } q_{i}=w^{j}(h-K|_{\Sigma})_{ij}.
\end{equation*}
According to the dominant energy condition this expression shows
that $\overline{R}$ is almost nonnegative, with only a divergence
term standing in the way.  In fact, the extra degree of freedom
given by the scalar $\phi$ will be used to remove the problematic
divergence term in the next section.  Moreover we have shown in
our investigation of the Penrose Inequality [20], in analogy with
the theory developed by Schoen and Yau [16] for the classical Jang
equation, that regular solutions of the modified Jang equation
exist in spherical symmetry away from apparent horizons if we
choose
\begin{equation}
\phi=\rho_{,s}
\end{equation}
where
\begin{equation*}
\partial_{s}=\frac{\sqrt{1-v^{2}}}{\sqrt{g_{11}}}\partial_{r},\text{
}\text{ }\text{ }\text{ }v=\frac{\sqrt{\phi^{2}
g^{11}}f_{,r}}{\sqrt{1+\phi^{2} g^{11}f_{,r}}}.
\end{equation*}
Note that
\begin{equation*}
s=\int_{0}^{r}\frac{\sqrt{1-v^{2}}}{\sqrt{g_{11}}}=
\int_{0}^{r}\sqrt{g_{11}+\phi^{2} f_{,r}^{2}}
\end{equation*}
is the radial arclength parameter for the $\overline{g}$ metric.
In particular we have\medskip

\textbf{Theorem 2 ([20]).}  \textit{Let $(M,g,k)$ be a spherically
symmetric initial data set satisfying the dominant energy
condition $\mu\geq |J|$. If a ball $B_{r}$ centered at the origin
does not contain an apparent horizon, then there exists a regular
solution $f$ in $B_{r}$ of the modified Jang equation with the
scalar $\phi$ given by (2.3).}

\begin{center}
\textbf{3. Existence of Apparent Horizons}
\end{center}\setcounter{equation}{0}
\setcounter{section}{3}

  Here we shall give the proof of Theorem 1, which will proceed by
contradiction.  Assume that the ball $B_{r}$ does not contain an
apparent horizon.  Then Theorem 2 guarantees the existence of a
regular solution to the generalized Jang equation with $\phi$
given by (2.3). In particular we must have $v(0)=0$, and $-1<v<1$.
Therefore $\phi$ and $\rho_{,r}$ are strictly positive on $B_{r}$
since
\begin{equation*}
4\sqrt{g^{11}}\frac{\rho_{,r}}{\rho}=2H_{S_{r}}=\theta_{+}+\theta_{-}>0,
\end{equation*}
as $\theta_{+}>0$ and $\theta_{-}>0$ when $B_{r}$ contains no
horizons.  Let
\begin{equation*}
m(r)=\sqrt{\frac{A_{\overline{g}}(S_{r})}{16\pi}}
\left(1-\frac{1}{16\pi}\int_{S_{r}}H_{S_{r},\overline{g}}^{2}\right)
=\frac{1}{2}\rho(r)(1-\rho_{,s}^{2}(r))
\end{equation*}
denote the Geroch energy [21] of a sphere $S_{r}$ inside the Jang
surface $\Sigma$.  A well-known calculation shows that
\begin{equation*}
m_{,s}=\frac{1}{4}\rho_{,s}\rho^{2}\overline{R},
\end{equation*}
so that the formula (2.2) for $\overline{R}$ yields
\begin{eqnarray}
m(r)=m(r)-m(0)&=&\int_{0}^{r}m_{,s}ds\nonumber\\
&=&\int_{0}^{r}\frac{1}{4}\rho_{,s}\rho^{2}\overline{R}ds\\
&=&\frac{1}{16\pi}\int_{B_{r}}\rho_{,s}\overline{R}d\omega_{\overline{g}}\nonumber\\
&\geq&\int_{B_{r}}\rho_{,s}(\mu-J(w)-(8\pi\phi)^{-1}\mathrm{div}
_{\overline{g}}(\phi q))d\omega_{\overline{g}},\nonumber
\end{eqnarray}
where $d\omega_{\overline{g}}$ is the volume element on the Jang
surface $\Sigma$.  We may then apply the divergence theorem (as a
result of the special choice of $\phi$ given by (2.3)) and the
calculation (see [20])
\begin{equation*}
\phi\overline{g}(q,n_{\overline{g}})d\sigma _{\overline{g}}
=-2\frac{\rho_{,r}v}{\sqrt{g_{11}}}
\left(\sqrt{g^{11}}\frac{\rho_{,r}}{\rho}v-k_{b}\right)d\sigma_{g}
\end{equation*}
where $n_{\overline{g}}$ is the unit outer normal to $S_{r}$ in
the $\overline{g}$ metric and $d\sigma_{\overline{g}}$,
$d\sigma_{g}$ are area elements, to obtain
\begin{eqnarray}
m(r)&\geq& 4\pi\int_{0}^{r}\rho_{,s}(\mu-J(w))\rho^{2}ds
-\frac{1}{8\pi}\int_{S_{r}}\phi\overline{g}(q,n_{\overline{g}})d\sigma_{\overline{g}}\\
&\geq&\frac{4\pi}{3}\int_{0}^{r}(\rho^{3})_{,r}dr\min_{B_{r}}(\mu-J(w))
+\frac{1}{4\pi}\int_{S_{r}}\frac{\rho_{,r}v}{\sqrt{g_{11}}}
\left(\sqrt{g^{11}}\frac{\rho_{,r}}{\rho}v-k_{b}\right)d\sigma_{g}\nonumber\\
&=&\frac{4\pi}{3}\rho^{3}(r)\min_{B_{r}}(\mu-J(w))+\frac{\rho_{,r}v}{\sqrt{g_{11}}}
\left(\sqrt{g^{11}}\frac{\rho_{,r}}{\rho}v-k_{b}\right)\rho^{2}(r).\nonumber
\end{eqnarray}
However since
\begin{equation*}
m(r)=\frac{1}{2}\rho(r)-\frac{1}{2}\left(\frac{1-v^{2}}{g_{11}}\right)\rho_{,r}^{2}\rho(r),
\end{equation*}
it follows that
\begin{eqnarray}
\frac{1}{2}\rho(r)&\geq&\frac{4\pi}{3}\rho^{3}(r)\min_{B_{r}}(\mu-J(w))
+\frac{1}{2}(1+v^{2})g^{11}\rho_{,r}^{2}\rho(r)-\frac{\rho_{,r}}
{\sqrt{g_{11}}}k_{b}\rho^{2}v(r)\nonumber\\
&=&\frac{4\pi}{3}\rho^{3}(r)\min_{B_{r}}(\mu-J(w))\\
&
&+\frac{1}{2}\rho^{3}\left(g^{11}\frac{\rho_{,r}^{2}}{\rho^{2}}-k_{b}^{2}\right)
+\frac{1}{2}\rho^{3}\left(k_{b}-\sqrt{g^{11}}\frac{\rho_{,r}}{\rho}v\right)^{2}\nonumber\\
&\geq&\frac{4\pi}{3}\rho^{3}(r)\min_{B_{r}}(\mu-J(w))+
\frac{1}{8}\rho^{3}(r)(H_{S_{r}}^{2}-(\mathrm{Tr}_{S_{r}}k)^{2}).\nonumber
\end{eqnarray}
Lastly because $\rho_{,r}>0$ we have
\begin{equation*}
\rho^{2}(r)\geq\frac{\int_{0}^{r}\sqrt{g_{11}}\rho^{2}}{\int_{0}^{r}
\sqrt{g_{11}}}=\frac{1}{4\pi}\frac{\mathrm{Vol}(B_{r})}{\mathrm{Rad}(B_{r})},
\end{equation*}
and hence
\begin{equation*}
\frac{3}{2}\frac{\mathrm{Rad}(B_{r})}{\mathrm{Vol}(B_{r})}\geq\min_{B_{r}}(\mu-J(w))
+\frac{3}{32\pi}\theta_{+}\theta_{-}(r).
\end{equation*}
We conclude that if (1.1) holds, then $B_{r}$ must contain an
apparent horizon.

\begin{center}
\textbf{4. Properties of the Misner-Sharp Energy}
\end{center}\setcounter{equation}{0}
\setcounter{section}{4}

  The Misner-Sharp energy [22] is widely regarded as the
correct measure of quasilocal energy contained in centered
spacelike 2-spheres in spherically symmetric spacetimes. When
evaluated on a sphere $S_{r}$ it takes the form
\begin{equation*}
E_{r}=\sqrt{\frac{A(S_{r})}{16\pi}}\left(1-\frac{1}{16\pi}\int_{S_{r}}\theta_{+}\theta_{-}\right),
\end{equation*}
which also happens to be the expression for the Hawking energy
[23] of a spacelike 2-surface in an arbitrary spacetime.  Here we
would merely like to point out that the arguments of the previous
section immediately imply positivity and monotonicity properties
for the Misner-Sharp energy.  To see this, let $B_{r_{1}r_{2}}$
denote the region between two concentric spheres $S_{r_{1}}$ and
$S_{r_{2}}$ with $r_{2}>r_{1}$.  We will refer to this region as
untrapped if $\theta_{+}\theta_{-}>0$ throughout. For definiteness
let us assume that both $\theta_{+}>0$ and $\theta_{-}>0$.  Then
$H_{S_{r_{1}}}\neq 0$ implies that
$|H_{S_{r_{1}}}^{-1}\mathrm{Tr}_{S_{r_{1}}}k|\leq 1$, so in
analogy with Theorem 2 there exists a regular solution of the
modified Jang equation with $\phi$ given by (2.3) and such that
$v(r_{1})=H_{S_{r_{1}}}^{-1}\mathrm{Tr}_{S_{r_{1}}}k$ (see [20]).
Note that this does not exclude the possibility that $S_{r_{1}}$
and/or $S_{r_{2}}$ are apparent horizons, and if this is the case
then we impose the restriction that they can be either future or
past but not both, which ensures that $H_{S_{r_{1}}}\neq 0$.
Therefore we may follow precisely the same arguments presented in
(3.1), (3.2), and (3.3) while keeping the two middle terms of
(2.2), to find that
\begin{eqnarray*}
E_{r_{2}}-E_{r_{1}}&=&\frac{1}{16\pi}\int_{B_{r_{1}r_{2}}}
\rho_{,s}(16\pi(\mu-J(w))+|h-K|_{\Sigma}|_{\overline{g}}^{2}+2|q|_{\overline{g}}^{2})d\omega_{\overline{g}}\\
&
&+\frac{1}{8}\rho^{3}(r_{2})(\mathrm{Tr}_{S_{r_{2}}}k-v(r_{2})H_{S_{r_{2}}})^{2}
-\frac{1}{8}\rho^{3}(r_{1})(\mathrm{Tr}_{S_{r_{1}}}k-v(r_{1})H_{S_{r_{1}}})^{2}.
\end{eqnarray*}
But since $v(r_{1})=H_{S_{r_{1}}}^{-1}\mathrm{Tr}_{S_{r_{1}}}k$ we
obtain
\begin{equation*}
E_{r_{2}}\geq E_{r_{1}}.
\end{equation*}
Conversely, if both $\theta_{+}<0$ and $\theta_{-}<0$ then the
same arguments with
$v(r_{2})=H_{S_{r_{2}}}^{-1}\mathrm{Tr}_{S_{r_{2}}}k$ give
\begin{equation*}
E_{r_{2}}\leq E_{r_{1}}.
\end{equation*}
We have thus found\medskip

\textbf{Theorem 3.}  \textit{Let $(M,g,k)$ be a spherically
symmetric initial data set satisfying the dominant energy
condition $\mu\geq |J|$.  Then the Misner-Sharp energy is always
monotonic on untrapped regions.  In particular, the Misner-Sharp
energy of a centered 2-sphere not enclosing any apparent horizon
is nonnegative $E_{r}\geq 0$, and the Misner-Sharp energy of a
centered two-sphere enclosing the outermost apparent horizon
$S_{r_{0}}$ satisfies the lower bound
$E_{r}\geq\sqrt{A(S_{r_{0}})/16\pi}$. Furthermore if $E_{r}=0$ or
$E_{r}=\sqrt{A(S_{r_{0}})/16\pi}$, then $(B_{r},g,k)$
(respectively $(B_{r_{0}r},g,k)$) arises from a spacelike
hypersurface in the Minkowski (respectively Schwarzschild)
spacetime.}\medskip

  These observations concerning the Misner-Sharp energy have
previously been established by Hayward in [24] (see also [25])
using different methods, although the rigidity result appears to
be new (for details see [20]). The novelty of our method lies with
the fact that it can easily be generalized to the nonspherically
symmetric case, once a general existence theory for the modified
Jang equation has been obtained.  When this is done, an expanded
version of Theorem 3 would give new positivity and monotonicity
properties for the Hawking energy, and would lead to a proof of
the Penrose Inequality [18] for general initial data.

\begin{center}
\textbf{References}
\end{center}

\noindent[1]\hspace{.06in}  K. Thorne, \textit{Magic without
Magic: John Archibald Wheeler}, edited by J.
Klauder\par\hspace{.06in}(Freeman, San Francisco, 1972), pp.
231.\medskip

\noindent[2]\hspace{.06in}  S. Hawking, G. Ellis, \textit{The
Large Scale Structure of Spacetime} (Cambridge
Univ.\par\hspace{.06in}Press, Cambridge, 1973).\medskip

\noindent[3]\hspace{.06in}  W. Israel, Phys. Rev. Lett
\textbf{56}, 789 (1986).\medskip

\noindent[4]\hspace{.06in}  P. Bizon, E. Malec, N. O'Murchadha,
Phys. Rev. Lett. \textbf{61}, 1147 (1988).\medskip

\noindent[5]\hspace{.06in}  P. Bizon, E. Malec, N. O'Murchadha,
Class. Quantum Grav. \textbf{6}, 961 (1989).\medskip

\noindent[6]\hspace{.06in}  J. Guven, N. O'Murchadha, Phys. Rev. D
\textbf{56}, 7658 (1997).\medskip

\noindent[7]\hspace{.06in}  J. Guven, N. O'Murchadha, Phys. Rev. D
\textbf{56}, 7666 (1997).\medskip

\noindent[8]\hspace{.06in}  E. Malec, Phys. Rev. D \textbf{49},
6475 (1994).\medskip

\noindent[9]\hspace{.06in}  D. Eardley, J. Math. Phys.
\textbf{36}, 3004 (1995).\medskip

\noindent[10] E. Malec, N. O'Murchadha, Class. Quantum Grav.
\textbf{21}, 5777 (2004).\medskip

\noindent[11] T. Zannias, Phys. Rev. D \textbf{45}, 2998
(1992).\medskip

\noindent[12] T. Zannias, Phys Rev D \textbf{47}, 1448
(1993).\medskip

\noindent[13] R. Schoen, S.-T. Yau, Comm. Math. Phys. \textbf{90},
575 (1983).\medskip

\noindent[14] S.-T. Yau, Adv. Theor. Math. Phys. \textbf{5}, 755
(2001).\medskip

\noindent[15] R. Beig, N. O'Murchadha, Phys. Rev. Lett.
\textbf{66}, 2421 (1991).\medskip

\noindent[16] R. Schoen, S.-T. Yau, Comm. Math. Phys. \textbf{79},
231 (1981).\medskip

\noindent[17] E. Malec, N. O'Murchadha, Phys. Rev. D \textbf{50},
R6033 (1994).\medskip

\noindent[18] H. Bray, M. Khuri, \textit{PDE's which imply the
Penrose Conjecture},  preprint,
\par\hspace{.01in} arXiv:0905.2622.\medskip

\noindent[19] P. Jang, J. Math. Phys. \textbf{19}, 1152
(1978).\medskip

\noindent[20] H. Bray, M. Khuri, \textit{A Jang equation approach
to the Penrose Inequality}, preprint,
\par\hspace{.01in} arXiv:0910.4785.\medskip

\noindent[21] R. Geroch, Ann. N.Y. Acad. Sci. \textbf{224}, 108
(1973).\medskip

\noindent[22] C. Misner, D. Sharp, Phys. Rev. \textbf{136}, B571
(1964).\medskip

\noindent[23] S. Hawking, J. Math. Phys. \textbf{9}, 598
(1968).\medskip

\noindent[24] S. Hayward, Phys. Rev. D \textbf{53}, 1938
(1996).\medskip

\noindent[25] G. Burnett, Phys. Rev. D \textbf{48}, 5688 (1993).

\bigskip\bigskip\footnotesize

\noindent\textsc{Department of Mathematics, Stony Brook
University, Stony Brook, NY 11794}\par

\noindent\textit{E-mail address}: \verb"khuri@math.sunysb.edu"

\end{document}